- The total number of words of the manuscript: 6748

- The number of words of the abstract: 220

- The number of figures: 8

- The number of tables: 0

# Graphical Abstract

## Characterisation of neonatal cardiac dynamics using ordinal partition network

Laurita dos Santos,Débora C. Corrêa,David M. Walker,Moacir F. de Godoy,Elbert E. N. Macau,Michael Small

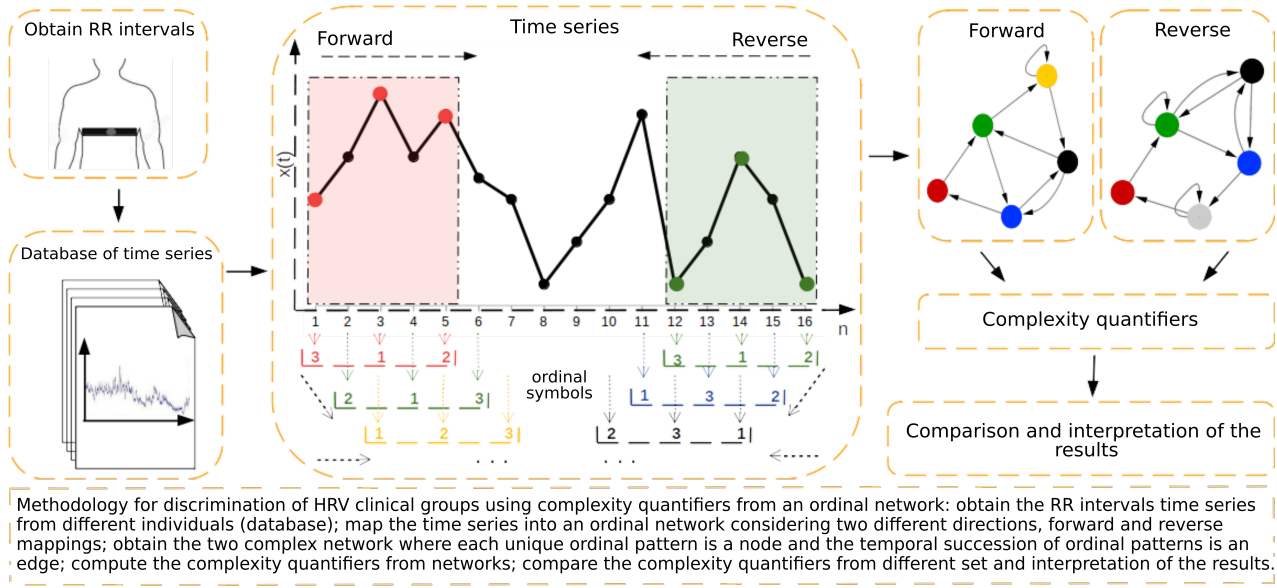

Methodology for discrimination of HRV clinical groups using complexity quantifiers from an ordinal network: obtain the RR intervals time series from different individuals (database); map the time series into an ordinal network considering two different directions, forward and reverse mappings; obtain the two complex network where each unique ordinal pattern is a node and the temporal succession of ordinal patterns is an edge; compute the complexity quantifiers from networks; compare the complexity quantifiers from different set and interpretation of the results.

# Characterisation of neonatal cardiac dynamics using ordinal partition network

Laurita dos Santos[a], Débora C. Corrêa[b,e], David M. Walker[b], Moacir F. de Godoy[c], Elbert E. N. Macau[d] and Michael Small[b]

[a]*Scientific and Technological Institute, Universidade Brasil, São Paulo, SP 08230-030, Brazil*

[b]*Complex Systems Group, Department of Mathematics and Statistics, The University of Western Australia, Crawley, WA 6009, Australia*

[c]*Department of Cardiology and Cardiovascular Surgery, Faculdade de Medicina de São José do Rio Preto, São José do Rio Preto, SP 15090-000, Brazil*

[d]*Institute of Science and Technology, Federal University of São Paulo, São José dos Campos, SP 12247-014, Brazil*

[e]*ARC Industrial Transformation Training Centre (Transforming Maintenance through Data Science), The University of Western Australia, Crawley, WA 6009, Australia*

## ABSTRACT

The maturation of the autonomic nervous system (ANS) starts in the gestation period and it is completed after birth in a variable time, reaching its peak in adulthood. However, the development of ANS maturation is not entirely understood in newborns. Clinically, the ANS condition is evaluated with monitoring of gestational age, Apgar score, heart rate, and by quantification of heart rate variability using linear methods. Few researchers have addressed this problem from the perspective nonlinear data analysis. This paper proposes a new data-driven methodology using nonlinear time series analysis, based on complex networks, to classify ANS conditions in newborns. We map 74 time series given by RR intervals from premature and full-term newborns to ordinal partition networks and use complexity quantifiers to discriminate the dynamical process present in both conditions. We obtain three complexity quantifiers (permutation, conditional and global node entropies) using network mappings from forward and reverse directions, and considering different time lags and embedding dimensions. The results indicate that time asymmetry is present in the data of both groups and the complexity quantifiers can differentiate the groups analysed. We show that the conditional and global node entropies are sensitive for detecting subtle differences between the neonates, particularly for small embedding dimensions ($m < 7$). This study reinforces the assessment of nonlinear techniques for RR intervals time series analysis.

## 1. Introduction

Heart rate variability (HRV) is an important non-invasive measurable physiological state related to homeostasis, which is defined as the self-regulating ability of a system or organism to maintain stability [1, 2]. HRV describes the variation of cardiac cycles represented by difference between two R waves (RR or NN intervals) in the electrocardiographic signal. This variation reflects contributions of the two branches of the autonomic nervous system (ANS): sympathetic and parasympathetic activities [3]. Particularly, the cardiac rhythm of the infants enable to elucidate the state of the ANS maturation, especially considering premature and full-term newborns [4, 5]. In terms of ANS maturation, distinguishing groups of time series from newborns (premature and full-term) is more difficult than distinguishing healthy newborns and healthy young adults, as there is a significant difference of age and degree of the maturation process of the organism.

The maturing process of the autonomic functions occurs in the gestation period, and premature newborns appear to have a delayed maturation of HRV when compared with full-term newborns [6]. In the gestation period, the foetal heart rate has neural and non-neural influences, such as parasympathetic and sympathetic innervations, a rise of parasympathetic influence, metabolic sensitivity, and other processes [7]. Thus, the measurement of HRV provides information about autonomic regulation of an individual, which reflects physiological fluctuations from two branches of inputs of the ANS. These oscillations have been related to nonlinear control systems representing the versatility of the healthy

lauritas9@gmail.com (L. dos Santos); debora.correa@uwa.edu.au (D.C. Corrêa); michael.small@uwa.edu.au (M. Small)
ORCID(s):

individual [8, 9, 10, 11]. This nonlinearity can be considered an intrinsic feature of the organism as a function of the state of the ANS maturation [12, 13].

A data-driven approach to classify RR intervals of premature newborns and full-term newborns can help to assess cardiovascular control and its specific features and limits. Furthermore, optimal heart control in newborns is not fully developed, and therefore future studies are required to elucidate its physiological phenomenon [12]. Today, the understanding of cardiac signals from newborns during the first month of their lives involves adjusted interpretation of electrocardiograms [14] and HRV analysis. Previous approaches attempting the automatic classification of RR intervals use both linear and nonlinear methods. Conventional linear methods include standard time domain and frequency domain approaches [1], which provide information about the periodic oscillations [4] of RR intervals. Time domain methods are the simplest way to determine the characteristics of successive RR intervals, such as the mean RR interval, the standard deviation RR interval and the difference between the longest and the shortest RR interval. These measures are useful to describe the instantaneous heart rate and cycle length [1]. Frequency domain methods are based on power spectral density (PSD), which uses spectral components, such as very low frequency, low frequency and high frequency, to provide information about autonomic modulations of the heart period [15, 16]. The main limitations of such techniques are imposed by the fluctuation of RR intervals time series, as they can provide similar or identical results for different profiles [1].

The nonlinear characteristics of physiological systems and how their interactions influence each other motivates nonlinear time series techniques. Nonlinear methods such as recurrence quantification analysis [17, 18], Poincaré plots [19], Shannon entropy quantification [20], correlation dimension [21], symbolic dynamics methods and time reversibility analysis [22], and physiology network [23] are good candidates for HRV analysis. All these methods allow associations between biological signals and their phenomenology in the organism, providing rich information about the systems' dynamics. Nonlinear methods differ from conventional linear HRV analysis because they assess other characteristics from the signal such as scaling, and nonlinear correlation properties [24]. Some advantages can be reported in the use of nonlinear approaches for HRV analysis, such as unpredictability, and complexity of the time series [25]. Here, we complement previous studies on this problem by proposing the application of nonlinear time series analysis from the perspective of complex networks. We focus, for the first time to the best of our knowledge, on the characterisation of time series from newborns (premature and full-term) based on the application of complex networks in biomedical data considering real-world clinical situations.

Representations of nonlinear systems using complex networks have attracted attention from researchers in different areas, such as human cardiac dynamics [26, 27, 28], epilepsy and seizure prediction [29, 30], climate models [31], finance [32], and investigation of determinism in time series [33]. In particular, complex networks provide a way to describe the evolution of system dynamics, where the entities describing a complex system are represented as nodes, and edges describe the intrinsic relationship among interconnected entities. These patterns of interactions can bring insights into the dynamical properties of the studied systems [34, 35]. In time series analysis, network representations try to capture the essence of the dynamics by representing the system state as nodes in the network and transitions between states as edges in the network [36]. For example, in a recurrence network, nodes are system states in the phase-space and are linked in the network if they are considered recurrent, that is, if they are part of the same defined neighbourhood in phase space [37]. Ordinal partition network methods also perform a partition in phase space by mapping the time series values to a sequence of symbols that define a ranking position of the amplitude values [34, 28]. Each segment will result in a permutation of symbols which are mapped to nodes in the network. Edges define the temporal sequence of such permutations that are characterising the system dynamics (defining a Markov transition network). As the length of the time series increases, new data tend to improve the estimation of the transitional probabilities instead of increasing the network size (as in the case of recurrence networks) [38].

The ordinal partition network therefore captures spatial dependencies of reconstructed phase space trajectories of the time series. Different dynamics of the time series will result in distinct exploration of phase space and consequently will result in different structural properties of the networks. Here, we propose to use the ordinal partition network and seek network properties that are sensitive to the structural changes of different nonlinear systems represented by RR intervals of premature newborns and full-term newborns. Measuring variability and complexity of such time series in their complex network representations is a novel and interesting way to characterise these complex systems.

We start by mapping the RR intervals time series from premature and full-term newborns into ordinal partition networks. Using this framework, we provide two contributions. First, we use the concept of surrogate data to show that complexity measures from the ordinal networks derived from RR intervals are helpful as discriminating statistics to inform the underlying dynamical process in the time series [39]. Second, we conduct experiments considering forward and reverse mappings of the RR time series [33] to evaluate the time irreversibility of these systems [40, 41]. In other words, the mapping of the time series into networks occurs both in forward and reverse directions. Our motivation is that the temporal asymmetry will provide a way to discriminate RR interval time series from the two classes in the database.

This paper is organised as follows: Section 2 describes the dataset and the strategy to map time series into ordinal networks and how we obtain the quantifiers reflecting the systems' complexity. Section 3 describes the viability of the ordinal network properties to display differences between the premature and full-term newborns time series; including the use of complexity quantifiers as a test for time irreversibility applied to the RR time series. Finally, we conclude the paper with a summary of our results.

## 2. Material and Methods

### 2.1. Dataset

The RR interval database comprises 74 time series or tachograms: 48 premature newborns (PNB) were an average $-27.4$ days of corrected age, and 26 full-term newborns (FNB) aged up to 8 days after their birth. The tachograms datasets were taken from previous studies of Transdisciplinar Nucleus for Chaos and Complexity Study (NUTECC/Brazil) [42]. Each tachogram was selected and classified (premature or full-term newborn) by a pediatric cardiologist, who was expert in cardiac signals analysis in infants and children. The dataset was approved by respective ethics committees and conducted in agreement with Helsinki's declaration for medical research involving human subjects, and all parents and/or guardians of the newborns gave their informed consent for participation in the study.

Data were collected by a Polar monitor (S800i or RS800) for approximately 15 minutes with a subject in a supine rest position. We compare the newborn's groups using $1,490$ RR intervals (cardiac cycles) for each individual. We pre-processed the data with an adaptive filter methodology to discard possible artifacts. The filtering algorithm first removes the RR intervals longer than $1,200$ ms and shorter than 350 ms as such intervals are inconsistent according to the physiology of the sinus rhythm. Based on the adaptive values of the average and standard deviation of the time series, RR intervals are replaced if they are more than 20% different than their adjacent points. The filtered time series that differ by more than 10% from unfiltered time series are not considered in the analysis. This adaptive filtering methodology was discussed in a previous study [43].

Figure 1 displays two examples of non-filtered RR intervals time series (Fig.1a from PNB group and Fig.1d from FNB group) and their corresponding filtered series (Figs.1b and 1e), respectively. The first non-filtered time series (Fig.1a) is an example with few artifacts, whereas the second one (Fig.1d) is an example with a large number of visible artifacts, which are not related to sinus rhythm. The filtering strategy clearly helps to reduce the artifacts. Figs.1c and Figs.1f shows the corresponding time series in more detail.

### 2.2. Ordinal partition network

The strategy to map the time series into a complex network is based on the ordinal partition presented in [28]. The authors showed that a univariate time series can be defined as a network with unique ordinal patterns (nodes), and edges linking nodes for which patterns occur sequentially in the observed data.

Given a tachogram $X = \{x_1, x_2, ..., x_N\}$, where $N$ is the number of the RR intervals, a series of embedding vectors $Z = \{z_1, z_2, ..., z_{N-(m-1)\tau}\}$ is obtained based on patterns of length $m$ and time lag $\tau$. We tested different values for $m$ and $\tau$ in our analysis. Each vector $z_i = \{x_i, x_{i+\tau}, x_{i+2\tau}, ..., x_{i+(m-1)\tau}\}$ is mapped to a symbol sequence according to the amplitude (y-axis) of each element in it. If two elements of $z_i$ have the same value ($x_i = x_j$), the symbol is assigned by order of appearance. A unique set of all symbol sequences $s$ is obtained in the time series, where $S = \{s_1, s_2, ..., s_{N-(m-1)\tau}\}$ for $s_i \in s$.

Each ordinal symbol in $s$ represents a node in the network $G = \{V, E\}$, where $V$ is the set of nodes ($|s|$) and $E$ is the set of edges. The adjacency matrix $A$ ($|V| \times |V|$) is the representation of the network, and $a_{i,j} > 0$ represents the link from node $i$ to node $j$. The weight of the directed edges in the network indicates the number of times that two adjacent symbols are presented sequentially in the time series.

In this work, the time series are mapped into a ordinal network considering two different directions, forward and reverse mapping, where the first one is the standard time sequence and the second one is the reversed time direction. It has been suggested that this technique is a reliable way to distinguish symmetrical from asymmetrical dynamical processes [10]. To illustrate the process, consider the time series with 16 RR intervals depicted in Figure 2. For $m = 3$ and $\tau = 2$, the first two forward ordinal patterns $s$ are $\{3, 1, 2\}$ and $\{2, 1, 3\}$, and the first two reverse ordinal patterns are $\{2, 1, 3\}$ and $\{2, 3, 1\}$ obtained by ranking the amplitude of the elements in the first two embedding vectors, respectively. Each unique ordinal pattern is a node of network and the temporal succession of ordinal patterns in the time series is an edge, as illustrated in 2e.

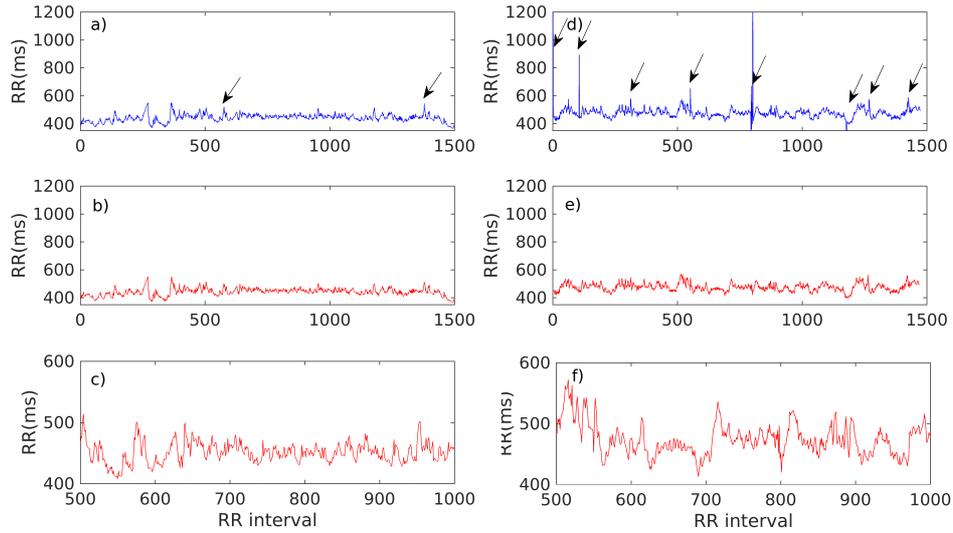

**Figure 1:** Examples of RR intervals time series before and after the adaptive filtering method: a) non-filtered time series of premature newborn with few visible artifacts (black arrows), b) filtered time series of premature newborn, c) detail of the time series shown in b) with 500 RR intervals, d) non-filtered time series of full-term newborn with a large number of artifacts (black arrows), e) filtered time series of full-term newborn, and f) detail of the time series shown in e) with 500 RR intervals.

Time-forward and corresponding time-reversed time series are statistically equally probable if they are time reversible [44]. In nonlinear time series analysis, time series reversibility analysis is used to identify systems for which dynamics look similar whether time runs forward or backward [45]. For instance, linear stochastic Gaussian processes (and their static nonlinear transformations) are statistically reversible, while systems of nonlinear dynamics are usually statistically irreversible [44]. There has been recent evidence that the statistical differences between networks constructed from forward and reverse mappings can be a strategy to quantify the time irreversibility [40, 33]. According to the literature, time irreversibility in cardiovascular signals is more frequent in short-time scales than in long-time scales [9, 10, 12]. These studies used time series with 256 cardiac cycles [9, 10] and 1000 cardiac cycles [12] to detect temporal asymmetries using indices such as Porta's index $P\%$, Guzik's index $G\%$ and Ehler's index $E$. Despite this interest, no one to the best of our knowledge has focused on the application of time irreversibility analysis to distinguish RR intervals from similar individuals' classes. Here, we use temporal asymmetry indexes from ordinal networks' properties to investigate individuals' RR interval classes and verify the feasibility of these properties to separate RR intervals from premature newborns and full-term newborns.

Once we get the networks derived from the forward and backward mapping of the time series, we propose to obtain quantifiers reflecting the complexity of the dynamics of the system from each adjacency matrix $A$: (i) permutation entropy ($h^{PE}$), (ii) conditional entropy ($h^{CPE}$), and (iii) global node entropy ($h^{GNE}$). These quantifiers can be used as a proxy of the complexity of the temporal structure in multiscale ordinal symbolic dynamics from time series, which can reveal changes on complexities of interbeat interval dynamics. We refer the reader to [28] for a comprehensive study regarding the use of these measures for ordinal partition networks when addressing properties of dynamical systems.

The permutation entropy ($h^{PE}$) is estimated by counting the probability of occurrence of each symbol in $S$ representing the state of the system:

$$h^{PE} = -\sum_{i}^{n} p_i \log p_i \tag{1}$$

where $n$ is the total of symbols, and $p_i$ is the relative frequency of each symbol in $S$. This measure describes the static complexity, randomness, or prediction accuracy of the system. In general, higher values of permutation entropy relate to lower prediction accuracies [46, 47]. The permutation entropy has been used in many biomedical applications, including epilepsy studies and cognitive neuroscience [48].

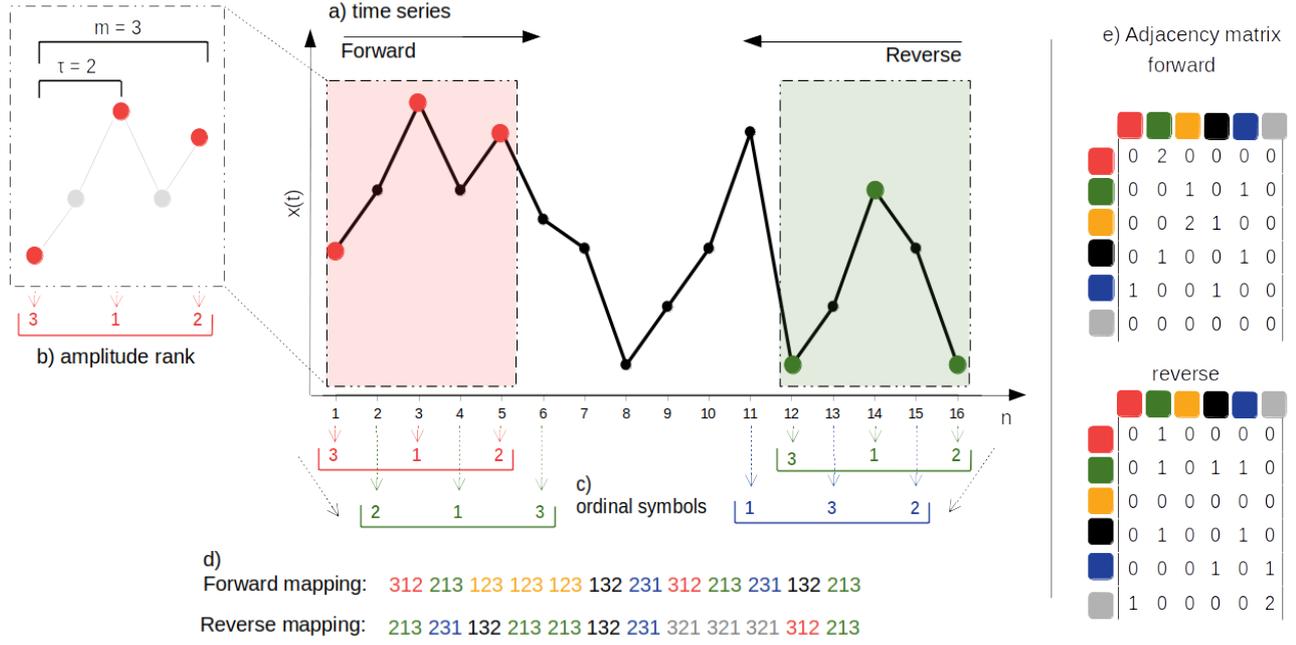

**Figure 2:** Example of the process for forward (left-to-right direction) and reverse (right-to-left direction) mappings of a time series in an ordinal network: a) find the first embedding vector with $\tau$ and $m$ according to the direction for mapping, b) determine its ordinal pattern by ranking the amplitudes of the elements in the window, c) obtain the ordinal symbols for all embedding vectors in both directions, d) obtain sequence of ordinal symbols based on both forward mapping and reverse mapping, e) compute the forward and reverse adjacency matrices using all ordinal patterns, where the edges represent the temporal sequence between nodes. Colours represent the unique ordinal symbols. Adapted from [28].

Conditional entropy ($h^{CPE}$) is an extension of permutation entropy [49], and quantifies the local uncertainty of each state in the model:

$$h^{CPE} = \sum_i \left( -p_i \sum_j p_{i,j} \log p_{i,j} \right) \tag{2}$$

where $p_i$ is the relative frequency of each symbol in $S$, and $p_{i,j}$ is an element of the stochastic matrix $P$ defined as the probability of a transition from symbol (state) $i$ to $j$:

$$p_{i,j} = \frac{a_{i,j}}{\sum_j a_{i,j}} \tag{3}$$

where $a_{i,j}$ are the elements of the adjacency matrix $A$.

For the calculation of both quantifiers in Equations 1 and 2, the adjacency matrix $A$ includes the possibility of self-loops in the network. Thereafter, a modified stochastic matrix $P$ (Eq. 3) is obtained by removing self-loops connections:

$$p_{i,j}^T = \begin{cases} 0 & \text{if } i = j \\ \dfrac{a_{i,j}}{\sum_{j, j \neq i} a_{i,j}} & \text{if } i \neq j. \end{cases} \tag{4}$$

The global node entropy ($h^{GNE}$) is related to fractal scaling exponent [28] and is defined as the transitional com-plexity value of $S$. According to [50], $h^{GNE}$ provides an average of the local transitional complexity of the attractor:

$$h^{GNE} = \sum_i p_i^* h_i^{LNE} \tag{5}$$

where $p_i^*$ is estimated by summing the weights for edges leaving node $i$ and dividing the result by the total sum of the edges weights of the ordinal network, and $h_i^{LNE}$ is the local node entropy defined as:

$$h_i^{LNE} = -\sum_j p_{i,j}^T \log p_{i,j}^T. \tag{6}$$

Statistical comparisons of $h^{PE}$, $h^{CPE}$, and $h^{GNE}$ obtained from a range of values for embedding dimension ($m$) and lag ($\tau$) regarding both forward and reverse mappings were carried out by Mann-Whitney test, with $p < 0.05$ for a 95% significance level. This statistical analysis was performed for intragroup comparisons and intergroup comparisons. In the first case, we use the network quantifiers for the same group's time series to investigate their temporal asymmetry. In the second case, we check whether these quantifiers can discriminate RR intervals from premature and full-term newborns. The code for mapping the time series into an ordinal partition network and its complexity quantifiers respectively is available at: https://github.com/laurita-santos/ordinal_partition_network.git.

## 2.3. Surrogate data

The surrogate framework was initially developed by Theiler and colleagues [39, 51]. The idea of this framework is to investigate explanations about the variability of a continuous time series. Suppose simple factors contributing to the variability can be discarded with some level of significance (for instance, variability is due to linear correlation or noise). In that case, one can justify a nonlinear analysis of the data. The basic approach consists of specifying a well-defined null hypothesis and generating an ensemble of surrogate data (new realisations of the time series) consistent with the null hypothesis. Then, one can check whether or not the null hypothesis can be rejected by looking at the distribution of a meaningful discriminating statistic estimated for the ensemble of surrogate data and the corresponding value for the original data. The surrogates framework has already been successfully used as a method to identify mechanisms of the system generating the observed data [51, 52, 53, 54, 55]. Here we will use the surrogate framework to check if variability measures in the time series help discriminate possibly nonlinear dynamics for the different groups in the dataset. The three null hypotheses proposed by Theiler [39, 51] are defined according to an increasing complexity assumption of the data: temporally independent data (Algorithm 0), linearly auto-correlated Gaussian noise (Algorithm 1), and static nonlinear transformation of linear Gaussian noise (Algorithm 2). For Algorithm 0 (Alg0), the null hypothesis states that the data is consistent with i.i.d. noise. In that case, surrogates can be random permutations of the original data, and the autocorrelation function at lags $\tau \geq 1$ will serve as a discriminating statistic. For Algorithm 1 (Alg1), the null hypothesis states that the data is consistent with linear filtered Gaussian noise. In this case, surrogates are obtained by shuffling the phases of the Fourier transform of the data. As rearranging the phases preserves the Fourier power spectrum, the surrogate data has the same autocorrelation characteristics (in addition to the probability distribution) as the original data. In that case, a nonlinear discriminating statistic is required, such as the mutual information at different lags. Algorithm 2 (Alg2) is an extension of Algorithm 1 for systems in which data is not normally distributed. Here, the null hypothesis states that data is consistent with a monotonic static nonlinear transformation of linearly filtered noise. In this case, the surrogate generated preserves the amplitude distribution and autocorrelation of the data. The most used procedure to generate surrogates consistent with this null hypothesis is Amplitude Adjusted Fourier Transform (AAFT) [39]. Mutual information at different lags, correlation dimensions and other nonlinear invariants can be used as the discriminating statistic. More recently, measures from complex networks have also been used as discriminating statistics for surrogate tests [56].

To say that the original data is statistically likely to be inconsistent with the null hypothesis, one must define an appropriate test statistic. One way to conduct a hypothesis test is to use a parametric criterion: the mean ($\mu_H$) and the standard deviation ($\sigma_H$) of a discriminating statistic computed from the surrogates are used to determine a significance level $\alpha = \frac{|Q_D - \mu_H|}{\sigma_H}$, where $Q_D$ is the value of the discriminant statistics for the original data [51]. However, for cases in which the statistical distribution is less likely to be Gaussian, the rank-order criterion is appropriate. This statistic-test ranks in increasing order the $N + 1$ discriminating statistics $Q_1, Q_2, ... Q_N$ computed for $N$ surrogates and $Q_D$ of the original data. If the original data is consistent with the null hypothesis, $Q_D$ has probability $1/(N + 1)$ to be the smallest or the highest value. Thus, the null hypothesis is rejected when $Q_D$ is the smallest or the higher value among all $(N + 1) Q$'s. This means that, for a two-sided hypothesis test at 95% significance level, a minimum of $N = 39$ surrogates are required.

In summary, the main steps of the proposed methodology are: a) map the time series into an ordinal network considering two different directions, forward and reverse mappings; b) compute the complexity quantifiers from adjacent matrix (networks); c) discriminate the complexity quantifiers of the RR intervals time series (PNB and FNB groups) from their surrogate data to assess the origin of variability in RR intervals; d) compare the complexity quantifiers obtained from forward and reverse mappings to characterise the system time irreversibility represented by the RR intervals time series.

## 3. Results

We present the results of our approach by first illustrating the viability of our methodology on chaotic time series from dynamical systems. Then we discriminate the variability in the RR intervals time series using the surrogate approach. Finally, we use complexity measures from original networks considering forward and reverse mapping to distinguish RR interval time series from the two groups.

### 3.1. Lorenz system and surrogate data

To show the use of surrogate data to unveil features about the variability of the time series, we generated time series from the standard Lorenz attractor under chaotic regime given by the following equations:

$$\begin{cases} \frac{dx}{dt} & = & \sigma(y - x) \\ \frac{dy}{dt} & = & x(\rho - z) - y \\ \frac{dz}{dt} & = & xy - \beta z \end{cases} \tag{7}$$

where $\sigma = 10$, $\rho = 28$, $\beta = 8/3$, and an incremental step 0.025 were used. For comparison, 10 time series were generated using a fourth-fifth Runge-Kutta method with randomised initial conditions, where each of $\{x(0), y(0), z(0)\} \in (0, 1)$. The peaks of $x-$component of Lorenz series (local maxima) with transients removed were sampled to form a time series with 1490 peaks (points). For the ordinal partition network, the parameter $m$ were selected according to the false nearest neighbor [57, 58]. For the peaks of $x-$component of Lorenz series the chosen parameters were $\tau = 1$ and $m = 3$.

For each of the three null hypotheses of surrogate data, a set of 100 surrogate series was constructed for each Lorenz series. As test statistics, we used the $h^{PE}$, $h^{CPE}$, and $h^{GNE}$ from the ordinal partition networks (only considering forward mapping in this case). Figure 3 shows the distribution of these three quantifiers for Lorenz series (peaks of $x-$component) and surrogate data. As expected, all hypotheses were rejected, the distribution of quantifiers for the data series differ from the distribution of quantifiers for the surrogate data, corroborating with the nonlinearity of the data. Clearly the $h^{PE}$, $h^{CPE}$, and $h^{GNE}$ for the surrogate time series show different distributions compared to the original time series. These results confirm previous statements that these quantifiers obtained from the complex network reflect the dynamics present in the analysed data [28].

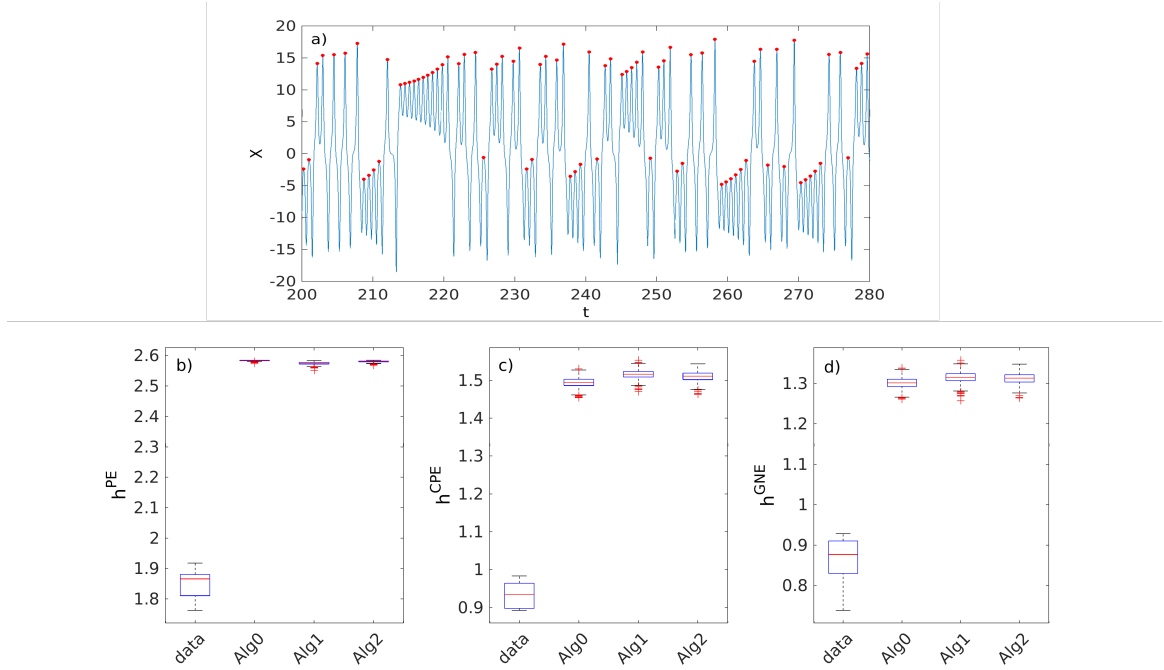

**Figure 3:** Distribution of quantifiers for time series given by the peak intervals of the Lorenz x-component and surrogate data: (a) time evolution of the $x-$component for Lorenz model (for $\sigma = 10$, $\rho = 28$, $\beta = 8/3$) and peaks marked with red dots detailed, (b) $h^{PE}$, (c) $h^{CPE}$, and (d) $h^{GNE}$ for peaks of the Lorenz time series and corresponding surrogate data, using parameters $m = 3$ and $\tau = 1$.

As expected, the test statistic distribution for the original data differs from the surrogate data distributions, showing the feasibility of this framework as hypothesis tests for nonlinear determinism.

### 3.2. Hypothesis tests for nonlinear determinism in RR interval time series

This work aims to assess the origin of variability in the RR intervals time series using the surrogate framework and the complexity quantifiers from ordinal partition network as test statistics. For each RR intervals time series we generate 100 surrogates time series consistent with null hypothesis given by Alg0 (non-correlated noise), Alg1 (linear correlated Gaussian noise), and Alg2 (nonlinear static transformation of linear Gaussian noise). To date, there is no established way to determine the optimum value for the embedding parameter $m$. It has been suggested that $m$ should be large enough to allow forbidden patterns, and small enough to promote accurate statistics of the transitions and minimise spurious results due to finite-size effects. Here, we map the RR time series into ordinal partition networks using parameter values $m = 1, 2, \ldots, 16$ and $\tau = 1, 2, 3, 4$ for both original and surrogate data, which cover common choices in the literature. Complexity quantifiers were obtained from the ordinal networks according to the description in the Section 2.2. We used a combination of the embedding $m$ and lag $\tau$ parameters to investigate the underlying dynamical process in the RR time series. We use the three quantifiers for each group of RR interval time series.

Figure 4 shows the statistical comparisons (p-values $p$) of the complexity quantifiers for different values of $m$ and $\tau$ considering RR interval time series of PNB and FNB groups and their corresponding Alg0 surrogate data. We observe that there is statistical difference ($p < 0.05$) for most of the parameters used, which supports the rejection of the null hypothesis for time series from both groups (left panel is the PNB data and right panel is the FNB data).

Comparisons between the original data and Alg1 surrogate data is presented in Figure 5. We note there is significant statistical difference for most of the parameters used between original data and Alg1 surrogate data for both premature group (left panel) and full-term group (right panel). The quantifier $h^{PE}$ of the PNB group (Figure 5a) detect differences between the original data and the Alg1 data for all values of parameters used, except $m = 7$. Interestingly, for the FNB group, we could accept the null hypothesis for other values of $m$. We observe that quantifiers $h^{CPE}$ and $h^{GNE}$ were consistent in distinguishing the original data from Alg1 surrogates when $m > 5$ up to $m < 12$ for all $\tau$ values for both two groups.

Similar behavior is observed when we compare original data to Alg2 surrogate data (Figure 6). There is an interval of parameter values (when $5 < m < 12$ for all $\tau$ used) in which the complexity quantifiers $h^{CPE}$ and $h^{GNE}$ reject the null hypothesis that variability of the data is consistent with nonlinear static transformation of linear Gaussian noise.

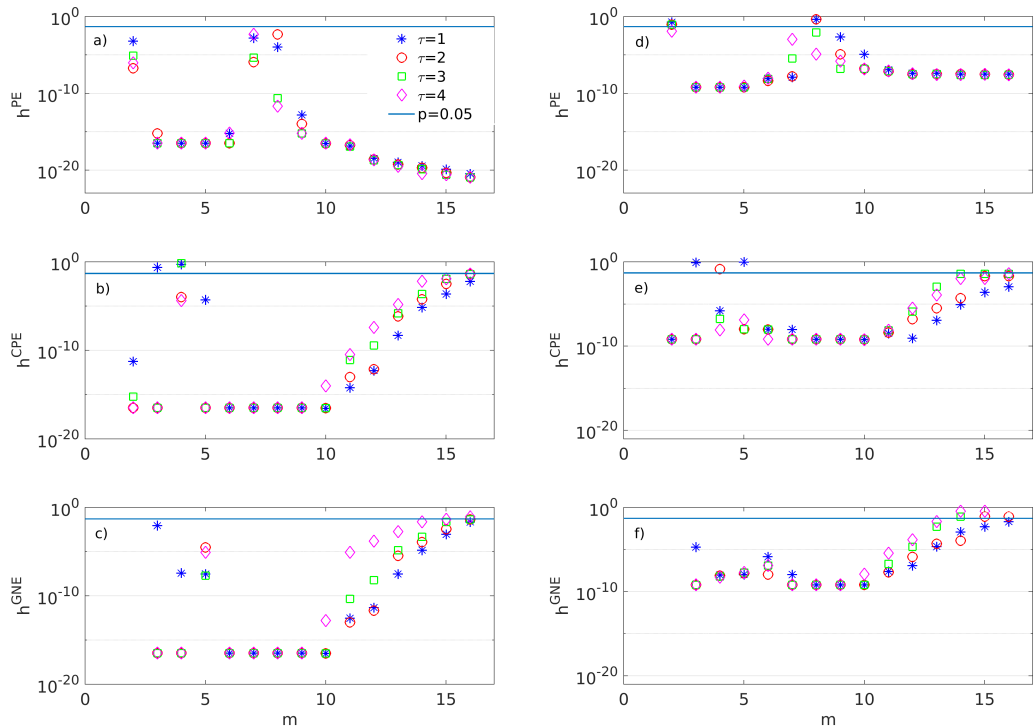

**Figure 4:** P-values results for the comparison between original data and Alg0 surrogate data considering the complexity measures $h^{PE}$, $h^{CPE}$ and $h^{GNE}$. Left panels show the results for the premature group and right panels show the results for the full-term group. Symbols below the blue line ($p = 0.05$) mean that there are statistical difference between the original data and the surrogate time series, supporting the rejection of the null hypothesis that variability of the data is consistent with i.i.d noise.

These results suggest the presence of nonlinearity for the system underlying the newborn time series, and that these quantifiers are sensitive and provide evidence regarding the difference between the dynamics present in the original data and the surrogates data, for a range of $m$ and $\tau$ used.

### 3.3. Comparison between forward and reverse mapping

The direction of the mapping (forward and reverse) of the complex network ordinal-based symbolic encoding is first used to assess the system time irreversibility represented by the RR intervals time series. It has been demonstrated that such forward and reverse encoding exhibited evidence for time irreversibility for low-order cycles of compression networks of chaotic time series [33]. Motivated by this study, here we investigate if forward and reverse mappings can provide information about time asymmetry of the distribution of the ordinal patterns found in the neonatal heart rate control. Here, each time series was mapped in two directions (forward and reverse), and we obtained the complexity quantifiers for each direction. We inspect evidence of time asymmetry using complexity quantifiers for intragroup time series (forward and reverse mappings for each time series of the same group), and intergroup time series (same network direction mapping for time series from PNB and FNB groups). A Mann-Whitney test is used to assess the significant difference in each case. This statistical test assumes a null hypothesis that there is no difference between forward and reverse mappings for each time series (intragroup comparison), and there is no difference between the time series from both premature and full-term newborns (intergroup comparison).

Figure 7 displays the intragroup statistical comparison between complexity measures obtained from forward and reverse network mapping for different values of $\tau$ and $m$. In the intragroup comparison for the premature newborns (PNB), we observed that there is statistical difference ($p < 0.05$) between both mappings mainly for the $h^{PE}$ and $h^{CPE}$ complexity quantifiers (see Fig.7a-b). For $h^{PE}$ and $h^{CPE}$ measures the null hypothesis is rejected when $\tau = 1$ and $\tau = 2$ for most of the $m$ values used. For the particular case $\tau = 1$, we verified that the $h^{PE}$ quantifier is able to detect the difference between the ordinal symbols obtained from the different mapping directions for all values of $m$.

For the $h^{CPE}$ quantifier, we observed a similar result for $\tau = 1$ (except for $m = 4$ and $m = 5$). These results suggest that these quantifiers reveal time asymmetry present in the RR intervals time series of the PNB group. For $\tau = 2$ we observed that the quantifiers $h^{PE}$ and $h^{CPE}$ reject the null hypothesis of symmetry at the level of 5% depending on the embedding value considered. For the $h^{GNE}$ the difference for $\tau = 1$ and $\tau = 2$ is also observed for larger values of $m$ ($m > 12$).

In the intragroup comparison for the full-term newborns (FNB), we observe that for $h^{PE}$ there is a statistically significant difference for $\tau = 1$ considering $m \le 8$ (Fig.7d), while $h^{CPE}$ suggest such difference for $m = 2, 4, 10$ up to 14 (Fig.7e). These results demonstrate that time asymmetry is also detected in the group of FNB, however for less combinations of values for $m$ and $\tau$ than for the group of prematures. The complexity quantifier $h^{GNE}$ displays a similar behavior for both intragroup comparisons: PNB (Fig.7c) and FNB (Fig.7f) groups, where the null hypothesis is accepted when $m < 11$ for all values of $\tau$, with exception for when $\tau = 1$ and $m = 3$ for the PNB group.

Figure 8 displays the intergroup statistical comparison between complexity measures obtained from forward (Fig.8a-c) and reverse mapping (Fig.8d-f). Here we observed that results for both forward and reverse mappings are similar (left and right panels), and they are consistent among the comparisons. The obtained p-values for the $h^{PE}$ measure demonstrate that there is statistically significant difference between PNB and FNB when $\tau = 1$ for different values of $m$, except $m = 3$ (forward mapping) and $m = 3$ and 4 (reverse mapping). Additionally, we observe in Figs.8a and 8d that the null hypothesis is rejected for all $\tau$ and $m$ values for $m \ge 9$. The $\tau$ is an important parameter to be properly established, particularly for $m < 9$.

Note that the quantifier $h^{CPE}$ shows similar results in Figs.8b and 8e for forward and reverse mapping, respectively. For both comparisons when $m \le 5$ the null hypothesis is rejected ($p < 0.05$), except for $\tau = 1$ and $m = 3$. For the quantifier $h^{GNE}$ (Figs. 8c and 8f) we observe that there is statistically significant difference for $m \le 6$ for both comparisons considering the values of $\tau$ used, except $\tau = 1$ and $m = 3$ (for PNB group) and $\tau = 1$ and $m = 4$ (for FNB group).

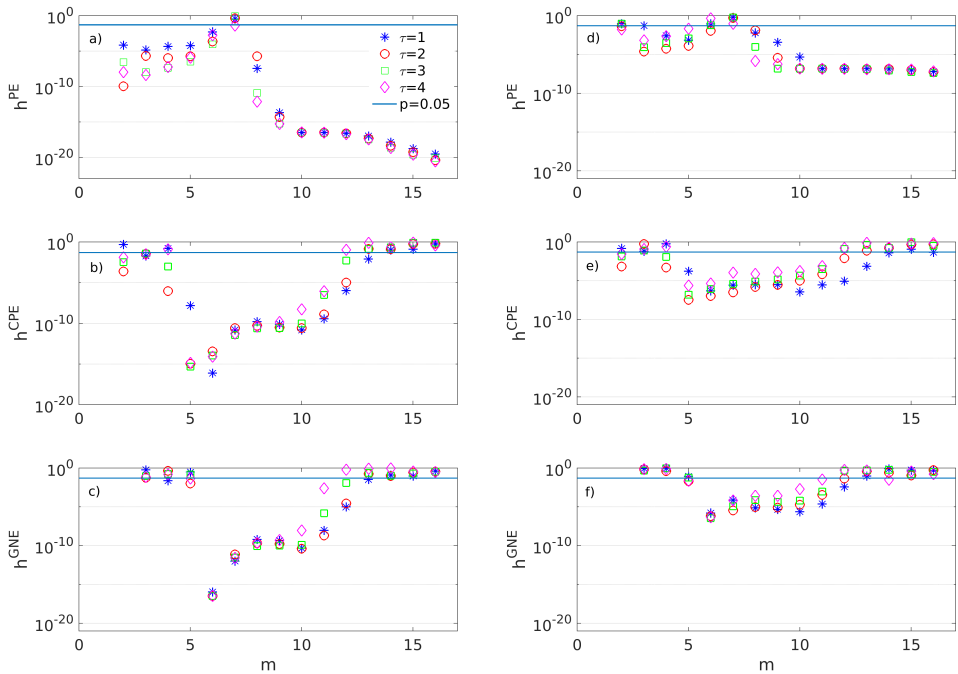

**Figure 5:** P-values results for the comparison between original data and Alg1 surrogate data considering the complexity measures $h^{PE}$, $h^{CPE}$ and $h^{GNE}$. Left panels show the results for the premature group and right panels show the results for the full-term group. Symbols below the blue line ($p = 0.05$) mean that there are significant statistical difference, then the null hypothesis is rejected. Symbols below the blue line ($p = 0.05$) mean that there are statistical difference between the original data and the surrogate time series, supporting the rejection of the null hypothesis that variability of the data is consistent with linear correlated Gaussian noise.

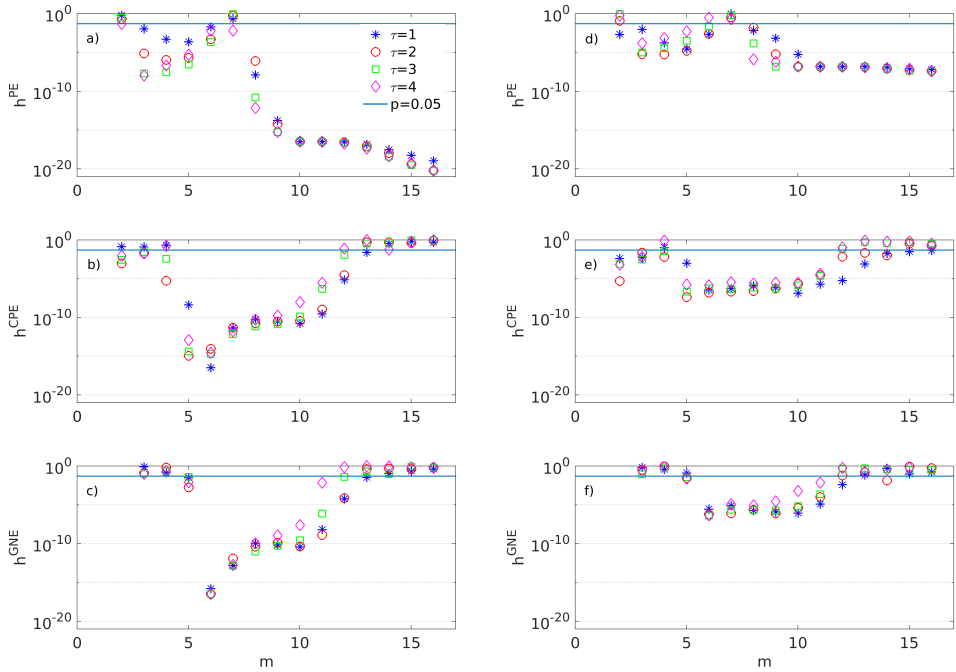

**Figure 6:** P-values results for the comparison between original data and Alg2 surrogate data considering the complexity measures $h^{PE}$, $h^{CPE}$ and $h^{GNE}$. Left panels show the results for the premature group and right panels show the results for the full-term group. Symbols below the blue line ($p = 0.05$) mean that there are statistical difference between the original data and the surrogate time series, supporting the rejection of the null hypothesis that variability of the data is consistent with nonlinear static transformation of linear Gaussian noise.

## 4. Discussion

Our results suggest that the time asymmetry is more evident in the PNB group in terms of numbers of combinations of $\tau$ and $m$ for which we can reject the null hypothesis. However, it is possible that such evidence for the presence of asymmetry in the PNB group is associated with the presence of accelerations and decelerations in the RR intervals time series in this group, which were related to sepsis [42, 59]. In our case, this clinical situation is a regular diagnosis for premature neonates under intensive care units. In this context, it is not possible to affirm that the time asymmetry can be associated only to the higher complexity of the PNB group.

In summary, our findings show that the complexity quantifiers indicate the time asymmetry in the RR interval time series for the PNB group and the FNB group, which is evidenced by the statistically significant differences shown in Fig. 7. This finding corroborates with the results in [12], in which they found that time asymmetry was present in heart rate oscillations in healthy eutrophic neonates. The authors also associated time asymmetry with nonlinear dynamics of the systems: the asymmetry seems to be higher of healthy newborns than healthy adults [12]. However, our results suggest that the time asymmetry is more evident in premature than in full-term newborns. The difference pointed out by the complexity quantifiers for the PNB group can be related to parasympathetic control's immaturity compared to the full-term newborns, which relative sympathetic dominance in neonatal heart rate control [11, 8].

Another finding of this study is that the complexity quantifiers detect differences in terms of maturation of autonomic functions between premature and full-term newborns. Some studies argue that HRV provides important information on the regulation of the cardiovascular system [18], and the ANS maturation for infants [42, 7, 4]. Our result shows that the complexity quantifiers based on ordinal partition network encoding display differences ($p < 0.05$) between the groups in the two directions of mapping the ordinal symbols (forward and reverse). We considered the complexity quantifiers are sensitive for detecting small differences between the neonates groups for embedding dimensions $m < 7$.

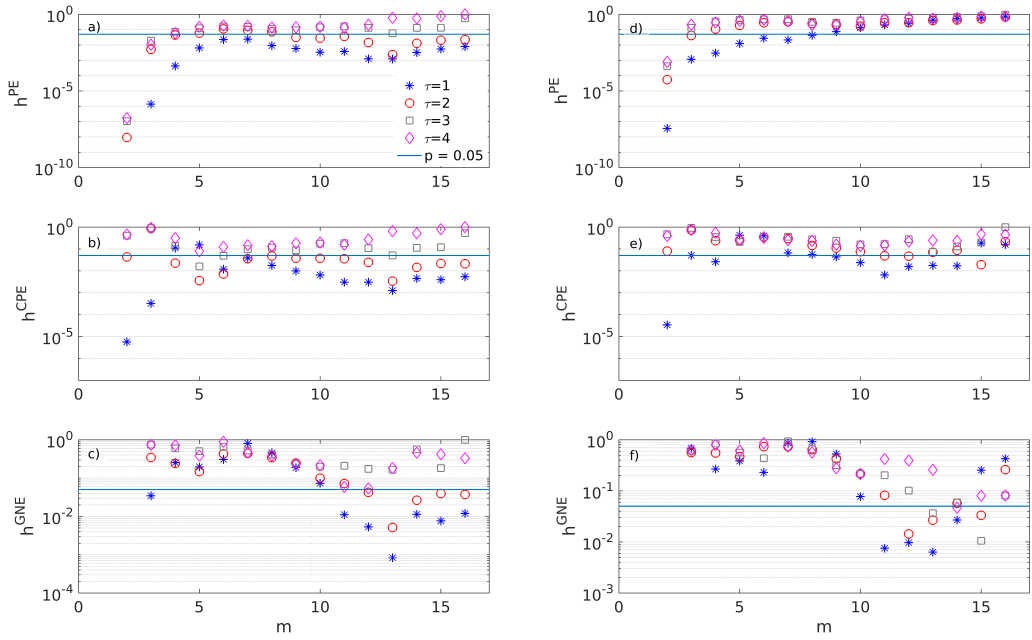

**Figure 7:** P-values results for the intragroup statistical comparison between forward and reverse mapping considering the complexity measures $h^{PE}$, $h^{CPE}$ and $h^{GNE}$. Left panels show the results for the premature group and right panels show the results for the full-term group. Symbols below the blue line ($p = 0.05$) mean that there are significant statistical difference, then the null hypothesis is rejected.

## 5. Conclusion

Some features about the underlying dynamical process in the HRV and how they impact on the physiological phenomenon are not completely understood. It is expected that levels of ANS maturation are different between premature and full-term newborn. This study looks into that question by assessing the complexity quantifiers obtained from ordinal partition network as classifiers for different ANS maturation periods.

Our results suggest that time asymmetry is present in the ordinal patterns of RR interval time series of both groups, as a proxy for the presence of nonlinear dynamics. In other words, the multiscale complexity measures obtained from forward and reverse mapping detect time asymmetry in the RR intervals time series. The results show that the time asymmetry is more present in the premature newborns than full-term newborns. This presence could be related to the accelerations and decelerations of RR intervals interpreted as sepsis. In terms of complexity quantifiers, we observed that they are sensitive to distinguish dynamics from the two groups, particularly considering $m < 7$, and this was consistent for forward and reverse mappings. These results reinforce that entropy-related measures from networks work well as complexity quantifiers, and can provide complementary information about the neonatal heart control from newborns.

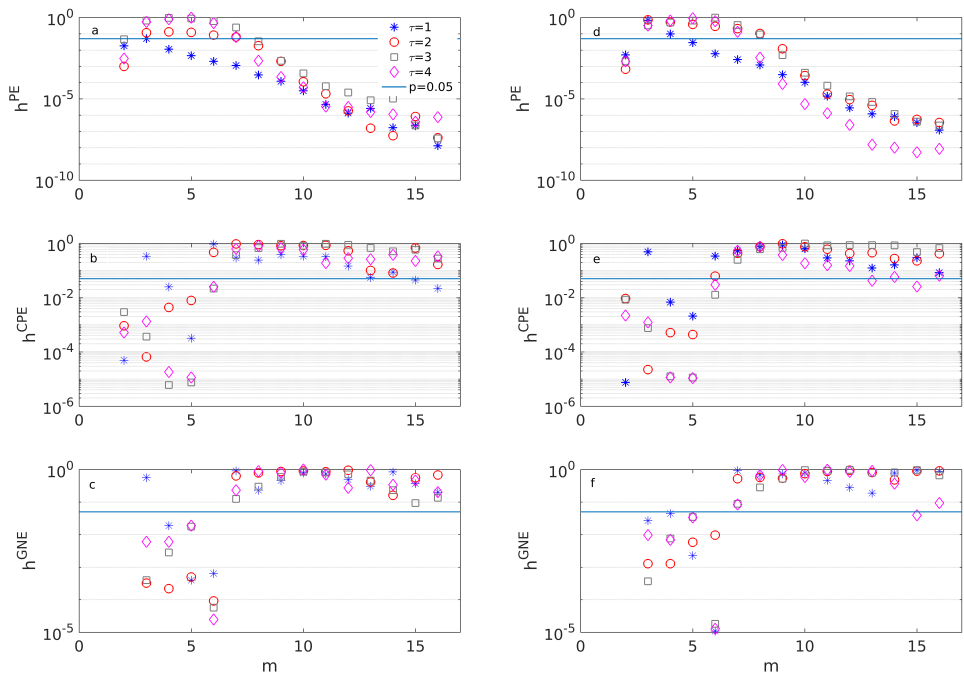

**Figure 8:** P-values results for the intergroup statistical comparison between PNB and FNB groups considering the complexity measures $h^{PE}$, $h^{CPE}$ and $h^{GNE}$. On the left panel we show the results for the forward mapping comparisons, and on the right panel for the reverse mapping comparisons. Symbols below the blue line ($p = 0.05$) mean that there are significant statistical difference, then the null hypothesis is rejected.

## Acknowledgements

L. dos Santos thanks the São Paulo Research Foundation - Fapesp (Grant no 2018/03517-8) and CNPq (Grant no 450133/2019-2) for their financial support. E.E. N.M. acknowledge the São Paulo Research Foundation - Fapesp (Grant no 2015/50122-0), CNPq and CAPES. M. Small, D. M. Walker, D. C. Correa (DP200102961) are partially supported by the Australian Research Council. D. C. Correa is also supported by the Australian Research Council through the Centre of Transforming Maintenance through Data Science (IC180100030).

## Conflict of interest statement

None declared.

---

**Laurita dos Santos** is Professor of Scientific and Technological Institute at Universidade Brasil. Her research interests are in the heart rate variability and nonlinear analysis applied to biomedical time series.

**Debora Correa** is a Research Fellow in the Department of Mathematics and Statistics at the University of Western Australia and one of the founding researchers of the ARC Industrial Transformation Training Centre for Transforming Maintenance through Data Science. Her research interests include nonlinear time series analysis and complex systems with applications in bioengineering and engineering systems.

**David Walker** is Senior Lecturer in the Department of Mathematics and Statistics at University of Western Australia. Research interests are in the areas of Complex Systems and Nonlinear Dynamics applied to biological systems.

**Moacir F. de Godoy** is Adjunct Professor of Cardiology. His main interests include Autonomic Nervous System, Heart Rate Variability and Complex Systems.

**Elbert E. N. Macau** is Full Professor at Universidade Federal de Sao Paulo. He serves as associate editor for International Journal of Bifurcation and Chaos, International Journal of Dynamics and Control, Frontiers in Computational Neural Sciences. His research interests include application of dynamical systems, data mind, time series analysis and complex networks.

**Michael Small** is Professor of Applied Mathematics at the University of Western Australia and CSIRO-UWA Chair of Complex Engineering Systems. His research interests include complex system, nonlinear dynamics, nonlinear time series analysis and chaos. He is Deputy Editor in Chief of the journal Chaos and Main Editor of Physica A.